\documentclass[12pt,leqno]{article}
\usepackage{graphicx}
\input amssym.tex
\begin{document}
\newtheorem{theo}{Theorem}[subsection]
\newtheorem{defi}{Definition}[subsection]
\newtheorem{cor}{Corollary}[subsection]
\newtheorem{lem}{Lemma}[subsection]
\newtheorem{pro}{Proposition}[subsection]
\newtheorem{que}{Exercise}[subsection]
\newtheorem{rem}{Remark}[subsection]
\newtheorem{exa}{Example}[subsection]
\def\cc{\c{c}}
\def\ii{\'\i}
\def\at{\~a}
\def\ot{\~o}
\newcommand{\nei}[2]{{#1_{(#2)}}}
\newcommand{\qed}{\hfill$\square$}
\def\pr{{\bf Proof: }}
\def\O{{\cal O}}
\def\M{{\cal M}}
\def\T{{\cal T}}
\def\N{{\Bbb N}}
\def\S{{\cal S}}
\def\Q{{\cal Q}}
\def\R{{\Bbb R}}
\def\v{\nu}
\def\C{\Bbb C}
\def\Qv{{\cal Q}_\v}
\def\nb{ \M/\M^2 }
\def\nbn{ \M'/\M'^2 }
\def\nbv{ \M^{\v}/\M^{\v+1} }
\title{Neighborhoods of Analytic Varieties in Complex Manifolds}
\author{C. Camacho, H. Movasati}
\maketitle
\tableofcontents
\section{Neighborhoods}
The systematic study of neighborhoods of analytic varieties was started by H. Grauert in his celebrated article \cite{gra}. In that article he considers a manifold $X$ and a negatively embedded submanifold $A\subset X$. He introduces $n$-neighborhood, $n\in {\Bbb N}$, of $A$ and studies when an isomorphism of two $n$-neighborhoods can be extended to an isomorphism of $(n+1)$-neighborhoods. He observes that obstructions to this extension problem lie in the first cohomology group of certain sheaves involving the normal bundle of $A$ in $X$. Using a version of Kodaira vanishing theorem (introduced by him in \cite{gra}) he shows that for a large $n$ these cohomology groups vanish and so he concludes that the germ of a negatively embedded manifold $A$ depends only on a finite neighborhood of it. These methods are generalized to a germ of an arbitrary negatively embedded divisor $A$ in \cite{hir} and \cite{lau}. In the case where $A$ is a Riemann surface embedded in a two dimensional manifold, by using Serre duality  we can say exactly which finite neighborhood of $A$ determines the embedding (see \cite{lau}). P. Griffiths in \cite{gri3} studies the problem of extension of analytic objects (fiber bundles, analytic maps, cohomology elements etc.) in $A$ to $X$. Again he introduces the finite extension of the object and he observes that obstructions lie in certain cohomology groups.

We must remark that  the above discussion leads  to formal extensions and isomorphisms. Grauert solves the convergency problem by geometrical methods. Laters, Artin's criterion (see \cite{art1}) on the existence of convergent solutions is used instead of Grauert's geometrical methods.

This is an expository article about negatively embedded varieties. The text is mainly based on Grauert's paper \cite{gra}, but we have used also 
the contributions of subsequent authors; for instance the Artin's criterion on the existence of convergent solutions is used instead of Grauert's geometrical
methods.
Our principal aim is to extend this study to the germ of foliated neighborhoods and singularities. In the first steps we will consider the most simple foliations which are transversal foliations. Next, foliations with tangencies and Poincar\' e type singularities will be considered.  


Let us fix the germ of an embedding $A\subset X$. We ask for the classification of germs of holomorphic foliations around $A$. Of course transversal foliations to $A$ will be of primary interest. We will try to generalize the Grauert step by step extension of isomorphisms to the case where the germ of embedding is foliated.

Artin in \cite{art1} after stating his extension and lifting theorems poses the following question:
Can one generalize these statements in various ways by requiring the map preserve extra structure, such as a stratification? In this paper this additional structure is a foliation. 
    
In the first section we will review some well-known facts about formal principle. Second section is devoted to Grauert's step by step extension of finite neighborhoods. The third section is devoted to the case of foliated  neighborhoods. The main result of this section is Theorem ~\ref{27jan01}. Also in this section we introduce the notion of formal equivalence of two foliated neighborhoods and prove Theorem ~\ref{26jan}. 

 For a mathematical autobiography of Grauert and also a brief history of complex analysis, the reader is referred to
 the interesting text \cite{rem1}.   
\subsection{Preliminary Notations and Definitions}
The basic notions of sheaf theory and cohomology theory, which can be found in 
\cite{gun3},\cite{grre},\cite{grre1}, are needed in understanding the following text. First, let us introduce some notations.

Let $A$ be a subvariety of a holomorphic variety $X$. We define:

\begin{itemize}
\item
$\O_X$, the structural sheaf of $X$, i.e. the sheaf of holomorphic functions on $X$;
\item
$\M$, the subsheaf of $\O_X$ consisting of elements that vanish at $A$;  
\item
$\nei{A}{*}=\O_X\mid_{A}$, $\nei{A}{*}$ is called the neighborhood sheaf of $A$;
\item
$A_{(\v)}= \O_X/\M^{\v}\mid_A$, $A_{(\v)}$ is called the $\v$-neighborhood of $A$;
\item
$\Qv= \nbv\mid_A$;
\item 
$\nei{\M}{\v}=\M/\M^{\v}\mid_A$.
\end{itemize}
Throughout the text, for a given sheaf $\S$ over $X$, when we 
write $x\in \S$ we mean that $x$ is a section of $\S$ in some open neighborhood in $X$ or it is an element in a stalk of $\S$ over $X$; being clear from the text which we mean.    
The natural inclusions 
\[
\cdots\subset\M^{\v+1}\subset\M^{\v}\subset\M^{\v-1}\subset\cdots\subset
\M
\]
give us the natural chain of canonical functions:
\[
\cdots \stackrel{\pi}{\rightarrow} A_{(\v+1)}\stackrel{\pi}{\rightarrow} A_{(\v)}\stackrel{\pi}{\rightarrow} A_{(\v-1)}
\stackrel{\pi}{\rightarrow}\cdots\stackrel{\pi}{\rightarrow} A_{(1)}
\]
We define
\[
A_{(\infty)}={lim}_{\infty\leftarrow \v} A_{(\v)}
\] 
In other words, every element of $A_{(\infty)}$ is given by the sequence 
\[
\ldots, f_{\v+1},f_{\v},f_{\v-1},\ldots,f_1\ \ \ f_v\in A_{(\v)} 
\]
\[
\pi(f_{\v+1})=f_{\v}
\]
The $\C$-algebra structure of $\nei{A}{\infty}$ is defined naturally.
$A_{(\infty)}$ is called the formal neighborhood of $A$ or the formal completion of $X$ along $A$.
There exists a natural canonical homomorphism
\[
\nei{A}{*}\rightarrow \nei{A}{\v}
\]
which extends to the inclusion 
\[
\nei{A}{*}\hookrightarrow \nei{A}{\infty}
\]
Define in the set 
\[
\tilde{\N}=\{1,2,3,\cdots,\infty,*\}
\]
the order
\[
1<2<3<\cdots<\infty<*
\]
we conclude that for any pair $\mu,\v\in \tilde{\N}, \mu\leq \v$ there exists a natural homomorphism
\[
\pi:\nei{A}{\v}\rightarrow\nei{A}{\mu}
\]
If no confusion is possible, we will not use any symbol for homomorphisms considered above.
Let us analyze the global sections of the above sheaves. Every global section of $\nei{A}{*}$ is a holomorphic function in a neighborhood of $A$. Let 
$g$ be a global section of $\nei{A}{\v},\v<\infty$. We can choose a collection of local charts $\{U_\alpha\}_{\alpha\in I}$ in $X$ covering $A$ and holomorphic functions $g_\alpha$ in $U_\alpha$ such that $g=g_\alpha$ in the sheaf $\nei{A}{\v}$. This means that 
\[
g_\alpha-g_\beta\in \M^\v\mid_{U_\alpha\cap U_\beta},\ \alpha,\beta\in I
\]
Conversely, every collection of $\{g_\alpha\}_{\alpha\in I}$ satisfying the above conditions defines a global section of $\nei{A}{\v}$.

We are interested in the germ of embedding $A\hookrightarrow X$. If 
$A\hookrightarrow X'$ is another embedding we denote by $A'$ its image. We fix the natural biholomorphism 
\[
\phi:A\rightarrow A'
\]
 obtained by the embeddings
$A\hookrightarrow X$ and $A\hookrightarrow X'$. This induces an isomorphism $\nei{\phi}{1}$ which is fixed from now on.
We always assume that the pairs $(X,A)$ and $(X', A')$ have the same local structure, i.e. for any $a\in A$ and its corresponding $a'=\phi(a)\in A'$ there is a local biholomorphism
\[
(X,A,a)\rightarrow (X',A',a')
\]
Notations related to $A'$ will be written by adding $'$ to the notations of 
$A$.  
\begin{defi}\rm
\label{com}
Let $\mu,\v\in \tilde{\N}, \mu\leq \v$.
The homomorphism $\nei{\phi}{\v}:\nei{A}{\v}\rightarrow \nei{A'}{\v}$ induces the homomorphism $\nei{\phi}{\mu}:\nei{A}{\mu}\rightarrow \nei{A'}{\mu}$, if the following diagram is commutative: 
\begin{equation}
\begin{array}{ccc}
\nei{A}{\v}& \stackrel{\nei{\phi}{\v}}{\rightarrow}  & \nei{A'}{\v} \\
\downarrow &  & \downarrow \\
\nei{A}{\mu} & \stackrel{\nei{\phi}{\mu}}{\rightarrow}  & \nei{A'}{\mu} \\
\end{array}
\end{equation}
We also say that $\nei{A}{\v}\rightarrow \nei{A'}{\v}$ extends $\nei{A}{\mu}\rightarrow \nei{A'}{\mu}$.  
\end{defi}
$\Q_1$ is the set of nilpotent elements of $\nei{A}{2}$ and so every 
homomorphism $\nei{\phi}{2}:\nei{A}{2}\rightarrow\nei{A'}{2}$ induces a homomorphism 
$\nei{\phi}{}:\Q_1\rightarrow\Q_1$. We also say that $\nei{\phi}{2}$ extends $\nei{\phi}{}$. 
\begin{defi}\rm
The homomorphism $\nei{\phi}{\infty}:\nei{A}{\infty}\rightarrow \nei{A'}{\infty}$ is called convergent if it takes $\nei{A}{*}$ into $\nei{A'}{*}$.
\end{defi}
{\bf Remark:}
In what follows, every homomorphism $\nei{A}{\v}\rightarrow\nei{A'}{\v}$ which we consider will be an extension of the fixed isomorphism  $\nei{A}{1}\rightarrow\nei{A'}{1}$ (Note that $\nei{A}{1}$ is the structural sheaf of $A$).

Let us analyze the geometrical interpretations of homomorphisms $\nei{\phi}{\v}:\nei{A}{\v}\rightarrow\nei{A'}{\v},\ \v\in \tilde{\N}$. The homomorphism 
$\nei{\phi}{\v}$ induces the identity $\nei{A}{1}\rightarrow\nei{A'}{1}$, and so, it takes the sheaf $\nei{\M}{\v}$ into $\nei{\M'}{\v}$. 

The following proposition gives us the 
local informations for analyzing the homomorphisms $\nei{\phi}{\v}:\nei{A}{\v}\rightarrow\nei{A'}{\v},\ \v\in \tilde{\N}$.
 \begin{pro}
\label{sara}
 Let $a\in A$ and $U$ be a small neighborhood of $a$ in $A$. Let also $a'=\phi(a)$ and $U'=\phi(U)$.
The following is true:
\begin{enumerate}
\item
 Every homomorphism (isomorphism) $\nei{\phi}{*}:\nei{A}{*}\mid_U\rightarrow \nei{A'}{*}\mid_{U'}$ which induces an isomorphism $\nei{A}{1}\mid_U\rightarrow \nei{A'}{1}\mid_{U'}$ is induced by a unique holomorphic (biholomorphic) map $(X',A',a')\rightarrow(X,A,a)$;
\item
Every homomorphism $\nei{\phi}{*}:\nei{A}{*}\mid_U\rightarrow \nei{A'}{*}\mid_{U'}$ which induces isomorphisms  $\nei{A}{2}\mid_U\rightarrow \nei{A'}{2}\mid_{U'}$, $\nei{A}{1}\mid_U\rightarrow \nei{A'}{1}\mid_{U'}$ is an isomorphism also;
\item 
Every homomorphism (isomorphism) $\nei{\phi}{\v}:\nei{A}{\v}\mid_U\rightarrow \nei{A'}{\v}\mid_{U'},\ \v>1$  is induced by a homomorphism (isomorphism) $\nei{A}{*}\mid_U\rightarrow \nei{A'}{*}\mid_{U'}$. 
\end{enumerate} 
\end{pro}
In the case where $a$ is a regular point of both $A$ and $X$, the proof of this proposition is easy. The proof in general uses simple properties of local
rings and their homomorphisms. The reader is referred to \cite{nag} for more informations about local ring theory. 

\pr
By Theorem 14, B of \cite{gun2} $\nei{\phi}{*}:\nei{A}{*}_a\rightarrow \nei{A'}{*}_{a'}$ is induced by a unique map $(X', a')\rightarrow(X,a)$. We must prove that this map takes $A'$ to $A$. Since $\nei{\phi}{*}$ induces an isomorphism  $\nei{A}{1}\mid_U\rightarrow \nei{A'}{1}\mid_{U'}$, it takes  the ideal of $A$ in $X$ to the ideal of $A'$ in $X'$. This implies that $(X',A',a')\rightarrow(X,A,a)$.

The second and third statements have a completely algebraic nature. To prove them
we use the following notations
\[
R=\nei{A}{*}_a\cong\nei{A'}{*}_{a'},\ I=\M_a,\cong M_{a'},\ \tau=\nei{\phi}{*},\ \tau_\v=\nei{\phi}{\v},\ \v\in \N
\]
( Note that $(X,A)$ and $(X',A')$ have the same local structure). Let us prove the second statement. 
Since $\tau_2:R/I^2\rightarrow R/I^2$ is an isomorphism and nilpotent set of $R/I^2$ is the set $I/I^2$, we have $I=\tau(I) +I^2$.  Let us prove that $\tau(I)=I$. Put
\[
R'=I/\tau(I)
\]
We have $IR'=R'$. Let $a_1,a_2,\ldots,a_r$ be a minimal set of generotors for $R'$. $a_r\in R'=IR'$ and so we have
\[
a_r=\sum_{i=1}^{r}s_ia_i,\ s_i\in I
\]
or
$(1-s_r)a_r$ lies in the ideal generated by $a_1,a_2,\ldots,a_{r-1}$. Since $1-s_r$ is inversible we get  a contradiction with this fact that no proper subset of $a_1,a_2,\ldots,a_r$ generates $R'$ (The used argument is similar to the proof of 
Nakayama's lemma (see \cite{gun2} A, Lemma 9)).

We have proved that $\tau(I)=I$. Since $\tau:R/I\rightarrow R/I$ is an isomorphism and $\tau(I)=I$, $\tau$ is surjective.
Now let us prove that $\tau$ is injective. Define
\[
R_n=\{x\in R\mid \tau^n(x)=0\}
\]
Since $R$ is a noetherian ring and $R_n\subset R_{n+1}$, there is a natural number $n_0$ such that $R_{n_0}=R_{n_0+1}=\cdots=R_*$. Now $\tau_*=\tau\mid_{R_*}$ is a surjective map from $R_*$ to $R_*$. By definition of $R_*$ we conclude that $R_*=0$ and so $R_1=0$. This means that $\tau$ is injective.  

Now let us prove the third statement.
 Let $x_1,x_1,\ldots,x_n$  form a basis for the vcecor space $\frac{\M_R}{\M_R^2}$, where $\M_R$ denotes the maximal ideal of $R$ (later in Proposition ~\ref{goftam} we will see that $x_i$'s form an embedding of $(X,a)$ in $(\C^n,0)$). 
 We can choose elements $f_1,f_2,\ldots,f_n$ in  $R$ such that $\tau_\v([x_i])=[f_i],\ i=1,2,\ldots,n$, where $[.]$ denotes the equivalence class. Now it is easy to verify that the homomorphism
\[
\tau: R\rightarrow R
\]
\[
f(x_1,x_2,\ldots,x_n)\rightarrow f(f_1,f_2,\ldots,f_n)
\]
induces the desired map. If $\tau_\v$ is an isomorphism then by the second part of the proposition $\tau$ is also an isomorphism.\qed

Now using Proposition ~\ref{sara} we can find geometrical interpretations of homomorphisms $\nei{A}{\mu}\rightarrow\nei{A'}{\mu},\ \mu\in \tilde{\N}$ as follows
\begin{enumerate}
\item
There exists an isomorphism
$\nei{\phi}{*}:\nei{A}{*}\rightarrow \nei{A'}{*}$ if and only if there exists a biholomorphism of some neighborhood of $A$ into some neighborhood of $A'$ in $X'$ extending $\phi:A\rightarrow A'$;
\item
Any isomorphism $\nei{\phi}{\v}:\nei{A}{\v}\rightarrow \nei{A'}{\v},1<\v\in \N$ is given by a collection of biholomorphisms $(U_\alpha, A)\rightarrow (U'_\alpha, A')$, where $\{U_\alpha\}_{\alpha\in I}$ ( resp. $\{U'_\alpha\}_{\alpha\in I}$) is an open covering of $A$ (resp. $A'$)
in $X$ (resp. $X'$), and such that $\phi_\alpha\circ\phi_\beta^{-1}$ is the identity up to holomorphic functions vanishing on $A$ of order $\v$;
\item
Any isomorphism $\nei{\phi}{\infty}:\nei{A}{\infty}\rightarrow 
\nei{A'}{\infty}$ is given by a collection of isomorphisms $\nei{\phi}{\v}:\nei{A}{\v}\rightarrow \nei{A'}{\v},\v\in\N$ such that for $\v \geq\mu$, $\nei{\phi}{\v}$ extends $\nei{\phi}{\mu}$.  
\end{enumerate}
The first statement justifies the name neighborhood sheaf adopted for $\nei{A}{*}$.

\subsection{Formal Principle}
The formal principle says: Every isomorphism from the formal
neighborhood $A\sb {(\infty)}$ of $A$ in $X$ to the formal neighborhood $A'\sb {(\infty)}$ of a subvariety $A'$ of a variety $X'$ implies the existence of a biholomorphism from an open neighborhood of $A$ in $X$ onto an open neighborhood of $A'$ in $X'$ (Note that we do not say that the formal isomorphism of neighborhoods is convergent). 
It is known that the formal principle does
not hold in every case.
V.I. Arnold in \cite{arn} has introduced a torus embedded in a complex manifold of dimension two with trivial normal bundle. The formal neighborhood of this 
torus is isomorphic with the formal neighborhood of the zero section of normal fiber bundle, but, there does not exist a biholomorphism between a neighborhood of the torus and of the zero section of normal fiber bundle.
However, the formal principle holds when the embedding of $A$ in $X$ 
has suitable properties of negativity \cite{gra}, 
or positivity \cite{hirs}. In \cite{art1} M. Artin proves the formal principle for singularities, i.e. $A=\{a\}$ is a single point  of a variety $X$. The next section is devoted to Artin's arguments.
\subsection{Formal Principle for Singularities}
This this section is devoted to the proof of formal principle for singularities following Artin's article \cite{art1}. When $A=\{a\}$ is a single point, following the literature we digress the notations as follows:
\begin{itemize}
\item
 $\O_{X,a}=\nei{A}{*}$, the local ring of $X$ at $a$; 
\item
$\hat{\O}_{X,a}=\nei{A}{\infty}$, the completion of the local ring $\O_{X,a}$;
\item
$\O_{X,a}^\v=\nei{A}{\v}$;
\item
For any local ring $R$, we denote by $\M_R$ its maximal ideal;
\item
$\C[[x]]$, the ring of formal series in $x$;
\item
$\C\{x\}$, the ring of convergent series in $x$.
\end{itemize}
In this section we need some basic informations about local rings. The reader can find these informations in \cite{nag}. 

Consider an arbitrary system of analytic equations
\begin{equation}
\label{30out00}
f_1(x,y)=0,f_2(x,y)=0,\ldots,f_k(x,y)=0
\end{equation}
where $f_1,f_2,\ldots,f_k$ are germs of holomorphic functions in $(\C^n\times\C^m,0)$.  
\begin{theo}
\label{artin}
(M. Artin \cite{art1}) Suppose that 
\[
\hat{y}(x)=(\hat{y}_1(x),\hat{y}_2(x),\ldots,\hat{y}_m(x))
\]
 are formal power series without constant term which solve ~(\ref{30out00}), i.e. 
\[
f(x,\hat{y}(x))= 0,\ f=(f_1,f_2,\ldots,f_k)
\]
 Let $c$ be a positive integer. There exists a convergent series solution
\[
y(x)=(y_1(x),y_2(x),\ldots,y_m(x))
\]
of ~(\ref{30out00}) such that 
\[
y(x)\equiv\hat{y}(x) \ \ modulo \ \M_{\C[[x]]}^c
\]
\end{theo}
Another way of stating the result is to say that the analytic solutions are dense in the space of formal solutions with its $\M_{\C[[x]]}$-adic metric (see \cite{nag} for definitions). 
\begin{cor}
\label{2nov00}
The formal principle holds for singularities, i.e. let $X$ and $X'$ be germs of holomorphic varieties at $a$ and $a'$ respectively. The isomorphism 
\[
\hat{\tau}: \hat{\O}_{X,a}\cong\hat{\O}_{X',a'}
\]
 of the formal completions  implies the isomorphism $X\cong X'$ of the germs of the varieties.  
\end{cor} 
The preceding corollary is the corollary 1.6 of \cite{art1}. It will be instructive to see how  the corollary can be obtained 
 from Theorem ~\ref{artin}.
\pr 
Without losing generality, we can assume that $X$ and $X'$ are germs of holomorphic varieties in $(\C^n,0)$ and $(\C^m,0)$, respectively. We use the following notations:
\[
\O_{X,0}=\C\{x\}/(f_1,\ldots,f_r)
\]
\[
\O_{X',0}=\C\{y\}/(g_1,\ldots,g_s)
\]
where $\C\{x\}$ is the ring of convergent series in $(\C^n,0)$. Let $\hat{p}_i(y)\in \C[[y]]$ represents the image $\hat{\tau}(x_i)$ of $x_i$ in $\C[[y]]$. The fact that $\hat{\tau}$ is a homomorphism implies that
\[
f_i(\hat{p}_1(y), \hat{p}_2(y),\ldots ,\hat{p}_n(y))\equiv 0\ \  \hbox{ modulo } (g_1,\ldots,g_s),\ i=1,2,\ldots,r
\]
i.e. there are formal series $\hat{\beta}_{ij}(y)\in\C[[y]]$ with
\[
f_i(\hat{p}_1(y), \hat{p}_2(y),\ldots ,\hat{p}_n(y))=\sum_{j=1}^{s}\hat{\beta}_{ij}(y)g_j(y)
\]
Now consider the system of holomorphic equations
\[
f_i(p_1,p_2,\ldots,p_n)-\sum_{j=1}^{s}{\beta}_{ij}g_j(y)=0,\ i=1,2,\ldots,r
\]
with unknown variables $y,p_i,\beta_{ij}$. This system has the formal solution 
$\hat{p}_i(y),\hat{\beta}_{ij}(y)$. Applying Theorem ~\ref{artin} with $c=2$, 
we obtain a homomorphism 
\[
\tau_1: \O_{X,0}\rightarrow \O_{X',0}
\]
which is congruent $\hat{\tau}$ modulo $\M_{\O_{X',0}}$. Let us prove that $\tau_1$ is an isomorphism.
With the same argument for $\hat{\tau}^{-1}$ we obtain $\tau_2$ congruent $\hat{\tau}^{-1}$ modulo $\M_{\O_{X,0}}$.  
Now $\tau_1\circ\tau_2:\O_{X,0}\rightarrow \O_{X,0}$ is congruent identity modulo $\M_{\O_{X',0}}$. The second part of Proposition ~\ref{sara} finishes the proof of our corollary.\qed
\subsection{Pseudoconvex Domains}
In this section we have introduced all preliminaries needed for the definition of negative embedding of varieties in the sense of Grauert. We start by the definition of plurisubharmonic functions and pseudoconvex domains $G\subset X$. 
\begin{defi}\rm
Let $G\subset \C^n$ be an open set  and $\phi:G\rightarrow \R$ a $C^2$-function. $\phi$ is called plurisubharmonic (resp. strongly plurisubharmonic) if
\[
\sum_{i,j=1}^{n}\frac{\partial^2\phi}{\partial z_i\partial \bar{z}_j}v_i\bar{v}_j\geq 0 \ \ (resp. > 0)\ \ \ \forall v=(v_1,v_2,\cdots,v_n)\not =0\in \C^n
\]
Let $X$ be an analytic space and $\phi:X\rightarrow\R$ a continuous function. $\phi$ is called plurisumharmonic (resp. strongly plurisubharmonic) if for any local chart $\psi:U\rightarrow V$, $U$ an open subset of $X$ and $V$ an analytic subset of an open subset of some $\C^n$, there exists a plurisubharmonic (resp. strongly plurisubharmonic) function $\check{\phi}$ on $V$ such that $\phi=\check{\phi}\circ\psi$.

The above definition is independent of the choice of local charts. 
\end{defi}
\begin{defi}\rm
Let $X$ be an analytic variety and $G$ a relatively compact open subset of $X$.  $G$ is called strongly pseudoconvex  if for every point $p$ in the boundary of $G$ there exist a neighborhood $U_{p}$ of $p$ and a real valued strongly plurisubharmonic $C^2$-function $\psi$ defined in $U_p$ such that
\[
G\cap U_p=\{x\in U\mid \psi(x)<0\}
\]
\end{defi}
For the definition of pseudoconvex spaces (= 1-convex spaces) the reader is referred to \cite{grre}. 
The next theorem shows that instead of local $C^2$-functions $\psi$, we can choose a global one.
\begin{theo}
\label{24jan01}
(\cite{gra} p. 338 Satz 2 \cite{grrepe} p. 228 Theorem 1.13) Let $G\subset X$ be a relatively compact strongly pseudoconvex domain. There exists a neighborhood $U$ of $\partial G$ and a  strongly plurisubharmonic
 $C^2$-function  $\psi$ in $U$ such that
\[
U\cap G=\{x\in U\mid \psi(x)<0\}
\]
\end{theo}
The idea is to use the partition of unity to patch together the local $\psi$'s and construct the global $\psi$. We need the above theorem to prove the following theorem:

\begin{theo}
(\cite{gra} p. 338 Satz 3)
Let $G\subset X$ be a relatively compact strongly pseudoconvex domain. Then there exists a compact set $K\subset G$ containing all nowhere discrete analytic compact subsets of $G$. 
\end{theo}
\pr
Let $\psi$ be as in Theorem ~\ref{24jan01} and 
\[
U_1=\{x\in U\mid -\epsilon <\psi(x)< 0 \}
\]
for a small $\epsilon$. We claim that $K=G-U_1$ is the desired compact set.
\\
Let $A$ be an analytic nowhere discrete compact subset of $G$ and $A\not \subset K$ or equivalently $A\cap U_1$ is not empty. $\phi$ has a maximum greater than $-\epsilon$ in $A\cap U_1$ which is a contradiction with this fact that $\psi$ is strongly pseudoharmonic.\qed

If $K$ is analytic, compact and nowhere discrete we say that $K$ is maximal. 
\begin{theo}
\label{dorahi}
(R. Narasimhan \cite{nar1} ). 
Let $G\subset X$ be a relatively compact strongly pseudoconvex domain. Then $G$ is holomorphically convex.
\end{theo}
Now the above theorem enables us to apply the Remmert reduction theorem to
strongly pseudoconvex domains. We recall some equivalent definitions of Stein
spaces.
\begin{theo}
The holomorphically convex space $X$ is Stein if and only if one of the following condition is satisfied:
\begin{enumerate}
\item
For any point $x\in X$ there exist holomorphic functions $f_1,f_2,\ldots,f_m$ on $X$ such that $x$ is an isolated point of the set $\{x\in X\mid f_1(x)=f_2(x)=\cdots=f_m(x)=0\}$;
\item
Holomorphic functions on $X$ separate the points of $X$, i.e. for any pair of points $x$ and $y$ in $X$ there exists a holomorphic function on $X$ such that $f(x)\not =f(y)$;
\item
$X$ does not contain nowhere discrete compact analytic subset;
\end{enumerate}
\end{theo}
The reader is referred to \cite{gun3} for the proof of the above theorem.
\begin{theo}
(Remmert reduction theorem)
Let $X$ be a holomorphically convex space. Then there exist a Stein space $Y$ and a proper surjective holomorphic map $\phi:X\rightarrow Y$ such that for any open set $U\subset Y$ and a holomorphic function $f$ in $\phi^{-1}(U)$ there exists a holomorphic function $g$ in $U$ such that $f=g\circ\phi$.
\end{theo}
{\bf Idea of the proof:} 
We define the relation $R$ in $X$ as follows:
\[
x_1 R x_2 \hbox{ if and only if } f(x_1)=f(x_2)\ \forall f\hbox{ holomorphic function in $X$}
\]
$R$ is an equivalence relation, and so, we can form the quotient space $Y=X/R$ and the continuous canonical map
\[
\phi: X\rightarrow Y
\]
We can also define a natural sheaf of holomorphic function $\O_Y$ on $Y$ as follows: For any open set $U\subset Y$,  $\O_Y(U)$ contains all holomorphic functions in $\phi^{-1}(U)$ which are constant on the fibers of $\phi$. We can easily verify that $\{\O_Y(U)\}_{U\subset Y}$ is a presheaf and so it defines a sheaf $\O_Y$. Now the principal point in the proof of Remmert reduction theorem is that:
\begin{itemize}
\item
$(Y,\O_Y)$ is a complex space.
\end{itemize}
$Y$ is holomorphically convex and global holomorphic functions separate the points of $Y$. Therefore $Y$ is Stein and the pair $(Y,\phi)$ is the desired pair.\qed

There are many contributions to complex analysis which are concerned with the following problem: When the quotient space of an equivalence relation in a complex space is again a complex space. The Grauert's direct image theorem plays an important role in these works. For more information the reader is referred to \cite{gra1} and its reference. The following statements are the consequences of Remmert reduction theorem:
\begin{enumerate}
\item
The pair $(\phi, Y)$ is unique up to biholomorphism;
\item
$\phi$ has connected fibers.
\end{enumerate}
The first is easy to prove. 
Let us check the second. Assume that $\phi^{-1}(x)$ is not connected and has two
connected components $A$ and $B$. In an open neighborhood of $\phi^{-1}(x)$ we can define a two valued function which takes $1$ in a neighborhood of $A$ and 0 in a neighborhood of $B$. This function is not a pullback of any holomorphic function in a neighborhood of $x$ in$Y$, which is a contradiction.

Before going to the next section let us state the most important sheaf theory property of strongly pseudoconvex domains:
\begin{theo}
(\cite{gra3})
\label{26jan01}
Let $G$ be a relatively compact strongly pseudoconvex domain in a complex manifold $X$ and $\S$ a coherent analytic sheaf on $G$. Then the cohomology groups $H^q(G, \S)$ are finite dimensional vector spaces for $q>0$.
\end{theo}
\subsection{Exceptional Sets}
Let $G\subset X$ be a relatively compact strongly pseudoconvex domain. By Theorem ~\ref{dorahi} $G$ is holomorphically convex, and so, we can apply Remmert reduction theorem to $G$ and obtain a Stein space $Y$ and a holomorphic map $\phi:G\rightarrow Y$. 
\begin{theo}
Let $G\subset X$ be a relatively compact strongly pseudoconvex domain and $\phi:G\rightarrow Y$ its Remmert reduction. Then the degeneracy set 
\[
A=\{x\in G\mid dim(\phi^{-1}(\phi(x))>0)\}
\]
is  the maximal compact analytic nowhere discrete subset of $G$.
\end{theo}  
\pr
 The subsets $\phi^{-1}\circ\phi(x),\in G$ are connected, and so by the definition, $A$ is nowhere discrete. We prove that $A$ is closed analytic set. The set $R=\{(x_1,x_2)\in X\times X\mid \phi(x_1)=\phi(x_2)\}$ is an analytic set and the projection on the first coordinate $\pi:R\rightarrow X$ is analytic. By \cite{gra1} Proposition 1 p.138 we know that
\[
\tilde{A}=\{x\in R\mid  dim(\pi^{-1}\pi(x))>0\}
\]
is a closed analytic set. Since $A=\pi(\tilde{A})$ and $\pi$ is proper, $A$ is also analytic closed set.

By Theorem \ref{24jan01}, there exists a compact set $K$ which contains all  compact analytic nowhere discrete subsets of $G$. For any $x\in A$, $\phi^{-1}\phi(x)$ is connected, and so by definition, is compact  nowhere  discrete subset of $A$. This implies that $\phi^{-1}\phi(x)\subset A$ and hence $A\subset K$. $A$ is closed set in the compact set $K$. Therefore $A$ is compact.\qed     

The Remmert reduction $\phi :G\rightarrow Y$ is proper and $A$ is compact so $\phi (A)$ is a compact analytic subset of $Y$. But $Y$ is Stein, and so, $\phi (A)$ is discrete set and $A$ is a union of compact connected analytic subsets $A_1,A_2,\ldots, A_r$ 
of $G$. In this case Remmert reduction substitute each $A_i$ with a point.
This leads us to the definition of exceptional sets.
\begin{defi}\rm
Let $X$ be an analytic variety and $A$ be a compact connected subvariety of $X$. $A$  is exceptional in $X$ if there exists an analytic variety $X'$ and a proper surjective holomorphic map $f:X\rightarrow X'$ such that
\begin{itemize}
\item
$\phi(A)=\{p\}$ is a single point;
\item
$\phi:X-A\rightarrow X'-\{p\}$ is an analytic isomorphism;
\item
For small neighborhoods $U'$ and $U$ of $p$ and $A$, respectively, $\O_{X'}(U')\rightarrow \O_X(U)$ is an isomorphism.
\end{itemize} 
\end{defi}
We also say that $A$ can be blown down to a point or 
is contractable or negatively embedded.
\begin{theo}
(Grauert,\cite{gra} Satz 5 p. 340)
Let $A$ be a compact connected analytic subset of $X$. $A$ is an exceptional set if and only if it has a strongly pseudoconvex neighborhood $G$ in $X$ such that $A$ is the maximal compact analytic subset of $G$.
\end{theo}
 \pr
Let us first suppose that $A$ is exceptional.  The analytic space $X'$ obtained by definition can be embedded in a $(\C^n,0)$ (definition of analytic sets). 
The neighborhood of $p$ in $X'$ given by
\[
U=\{x\in X'\mid z_1(x)\overline{z_1(x)}+\cdots +z_n(x)\overline{z_n(x)}<\epsilon\},\ \epsilon \hbox{ a small positive number }
\]
 is a pseudoconvex domain. Now it is easy to see that $G=\phi^{-1}(U)$ is the desired open neighborhood of $A$.

Now let us suppose that $A$ has a strongly pseudoconvex neighborhood $G$ in $X$ such that $A$ is the maximal compact analytic subset of $G$. Let $\phi:G\rightarrow X'$ be the Remmert reduction of $G$. We can see easily that $A$ is the degeneracy set of $\phi$ and $\phi(A)$ is a single point $p$. Since the fibers $\phi^{-1}\phi(x)$ are connected, $\phi$ is one to one map between 
$G-A$ and $X'-\{p\}$. Combining this and the property of $\phi$ in
 Remmert reduction theorem we can conclude that $\phi$ induces an isomorphism of stalks and so  it is a biholomorphism between  $G-A$ and $X'-\{p\}$. The third condition of an exceptional set can be read directly from Remmert reduction theorem.
\qed   

 When $A$ is a union of curves in a two dimensional manifold we have a
numerical criterion for contractablity of $A$.
\begin{theo}
\label{2nov2000}
Let $A$ be a compact connected one dimensional subvariety of a manifold $X$.
Suppose that $A$ contains only normal crossing singularities. Then $A$ is exceptional in $X$ if and only if the intersection matrix $S=[A_i.A_j]$ of $A$
in $X$ is negative definite, where $A=\cup A_i$ is the decomposition of $A$ into irreducible components. 
\end{theo}
This is Theorem 4.9 of \cite{lau}.

The concept of being exceptional is contained in which neighborhood of $A$? Let $A'$ be the image of another embedding $A\hookrightarrow X'$ of $A$. The following lemma gives an answer to this question. 
\begin{theo}
\label{vivoporela}
If $A$ is exceptional and there exists an isomorphism $\nei{\phi}{2}:\nei{A}{2}\rightarrow\nei{A'}{2}$  of 2-neighborhoods then $A'$ is also exceptional.  
\end{theo}
This is Theorem 6.12 of \cite{lau}, Satz 8 p.353 of \cite{gra} and Lemma 11 of \cite{hir}. The idea of the proof in the situation of Theorem ~\ref{vivoporela} is the following:

$\nb$ is the nilpotent subsheaf of $\nei{A}{2}$ and so every isomorphism of 2-neighborhoods induce an isomorphism of $\nb$'s. Therefore $A$ and  $A'$ have the same intersection matrix. Now  according to Theorem ~\ref{2nov2000}, $A'$ is also exceptional. \qed
\subsection{Formal Principle for Exceptional Sets}
For the first time formal principle was proved by Grauert in \cite{gra} for  compact  manifolds with codimension one. Extending his methods, Hironaka and Rossi in \cite{hir} proved the formal principal for exceptional sets.
\begin{theo}
(\cite{gra},\cite{hir})
\label{25jan01}
The formal principle holds for exceptional sets.
\end{theo}
To state a more general theorem we give the precise definition of modification.
\begin{defi}
A proper surjective holomorphic map $\phi:X\rightarrow Y$ of analytic varieties $X$ and $Y$ is called a modification if there are closed analytic sets $A\subset X$ and $Y\subset Y$ with codimension at least one such that
\begin{enumerate}
\item
$\phi(A)=B$;
\item
$\phi: X-A\rightarrow Y-B$ is biholomorphic;
\item
$A$ and $B$ are minimal with the properties 1 and 2.
\end{enumerate} 
\end{defi}
Note that for us an analytic space is always assumed to be reduced. A more
 general theorem about formal principle is the following:
\begin{theo}
(\cite{kos},\cite{anc})
\label{del2asir}
If $\phi:(X,A)\rightarrow (Y,B)$ is a modification with $A$ and $B$ compact then the formal principle holds for $(X,A)$ if and only if it holds for $(Y,B)$.   
\end{theo}  
Formal principle is true for singularities and so Theorem ~\ref{25jan01} is a consequence of Theorem ~\ref{del2asir}.
\subsection{Construction of Embedded Riemann Surfaces}
\newcommand{\modul}[1]{{\cal I}(#1)}
\newcommand{\proje}[1]{{\Bbb P}^{#1}}
 In this section we discuss various ways for constructing an embedding of a Riemann surfacce $A$ in a two dimensional manifold. The positive embeddings are abundant. They can be obtained by hyperplane sections of two dimensional algebraic manifolds. The first natural way to get a negative embedding is the following:

Let $A$ be a Riemann surface and $A\hookrightarrow X$ a positive embedding of $A$ in a two dimensional manifold, i.e. $A.A\geq 0$. Performing a blow up in a point $x$ of $A$ gives us another embedding of $A$ in a two dimensional manifold with self-intersection $A.A-1$. In fact, the new normal bundle of $A$ is $N.L_{-x}$, where $N$ is the normal bundle of $A$ in $X$ and $L_{-x}$ is the line bundle associated to the divisor $-x$.
Performing more blow ups in the points of $A$ gives us negative embeddings of $A$ with arbitrary self-intersection.

We have learned another way of changing the normal bundle of an embedding from P. Sad which goes as follows. The basic idea comes from \cite{cam2}.

Fix a germ of an embedding $(X,A)$ (for instance we can suppose $X=A\times \C$).
 Let $\S$ be the sheaf of local biholomorphisms $(X,A,x)\rightarrow (X,A,x),\ x\in A$ sending $A$ to $A$ identically. $\S$ is clearly a non-Abelian sheaf. We define an equivalence relation in $H^1(A, \S)$ as follows:  For $F=\{F_{ij}\},F'=\{F'_{ij}\}\in H^1(A, \S)$, 
 $F\sim F'$ if and only if there exists a collection of biholomorphisms $\{g_i\}$ such that
\[
F'_{ij}=g_i\circ F_{ij}\circ g_j^{-1}
\]
We define
\[
\modul{X}=H^1(A, \S)/\sim
\]
To each $F\in\modul{X}$ we can associate the line bundle $L_F=\{det( DF_{ij}\mid_{A})\}$.

Let $\{\psi_i\}$ be a collection of chart maps for the the germ $(X,A)$ and $F=\{F_{ij}\}\in H^1(A,\S)$. The new collection of transition functions
\[
\psi_i\circ F_{ij}\circ\psi_j^{-1}
\]
defines an embedding of $A$ with the normal bundle $L_FN$. We can see easily that two $F,F'\in
H^1(A,\S)$ give us the same embedding if and only if $F\sim F'$. Therefore we have
\begin{pro}
$\modul{X}$ is the moduli space of germs of all embeddings of $A$ in two dimensional manifolds. Moreover the line bundle of the embedding associated to $F\in\modul{X}$ is $L_F.N$.  
\end{pro}

Now let $X=N$ be a linear bundle. Consider another line bundle $M$ over $A$ with a meromorphic section $s$ of $M$. There is defined the biholomorphism
\[
\delta: N\rightarrow N.M
\]
\[
v\rightarrow v.s 
\]
which is well-defined out of the fibers passing through the zeros and poles of $s$. Now we can define
\[
\Delta: \modul{N}\rightarrow \modul{NM}
\]
\[
\{F_{ij}\}\rightarrow \{ \delta\circ F_{ij}\circ \delta^{-1}\}
\]
The line bundle associated to $F$ and $\Delta(F)$ are equal but the normal bundle of the embedding associated to $F$ is $L_FN$ and of $\Delta(F)$ is $L_FNM$. 

Another interesting method which can give us embedded Riemann surfaces is the action of groups.
Consider a subgroup $G$ of $Diff(\C^2,0)$ and denote by $G_0$ its linear group. After a blow up in $0\in \C^2$ we can consider $G$ as a group which acts in a neighborhood of $\proje{1}$, the projective line of blow up. Now we assume that $G_0$ is a Kleinian group which acts on $\proje{1}$. If $U_0$ is a region in $\proje{1}$ such that $A=U_0/G_0$ is a compact Riemann surface then it would be interesting to find a region $U$ in a neighborhood of $U_0$ in the blow up space such that $(U/G, A)$ is an embedding of $A$. For more information about Kleinian groups  the reader is referred to \cite{mas} and \cite{leh}.

\section{Obstructions to Formal Isomorphism}
\label{10feb01}
In this section we will identify the obstructions for the existence of an isomorphism between formal neighborhoods of $A$ and $A'$. We use the following trick: First we observe  when it is possible to get isomorphism of 2-neighborhoods of $A$ and $A'$ knowing the isomorphism of their normal bundles. Second when an isomorphism of $\v$-neighborhoods extends to an isomorphism of $(\v+1)$-neighborhoods. The applied methods are quite formal and can be found in \cite{gra},\cite{hir},\cite{lau}.     
\subsection{Notations and Exercises}
From now on, suppose that $A$ is an analytic subvariety of a manifold $X$. If $A$ and $X$ are smooth then  $N$ denotes the normal bundle of $A$ in $X$. The case where $A$ is a Riemann surface embedded in the complex surface $X$ will have our special attention. Let us start this section with some notations and exercises.
\begin{itemize}
\item
For any vector bundle $F\rightarrow X$, $F^*$ denotes its dual and $\underline{F}$
the sheaf of holomorphic sections of $F$;
\item 
For any analytic sheaf $\S$ on $X$, $\S^*$ denotes its dual i.e.,
\[
\S^*=\hbox{ the sheaf of $\O_X$-morphisms } \S\rightarrow \O_X
\]
\item
For any analytic sheaf $\S$ on $X$
\[
res(\S)=\S/\S. \M
\]
is called the structural restriction of $\S$ on $A$. Note that the sheaf theory restriction $\mid$ has nothing to do with the complex structure of $A$ but this restriction has. For instance the structural restriction of $\O_X$ to $A$ is $\O_A$. When there is no danger of confusion we will write the same symbol $\S$ instead of $res(\S)$; 

\item
For any analytic sheaf $\S$ on $X$
\[
\S(\v)=res(\S)\otimes_{\O_A} \Qv
\]
(Note that $\Qv$ is a $\O_A$-module sheaf); 
\item 
$\T$, the sheaf of holomorphic vector fields in $X$ (sections of the tangent bundle $TX$);
\item
 $\T_A$, the subsheaf of $\T$ consisting of vector fields tangent to $A$;
\item
When $A$ is a codimension one submanifold of $X$, by $N$ we denote the normal bundle of $A$ in $X$.
\end{itemize}
The reader will dominate the above notations if he succeed to verify the following simple facts:
\begin{itemize}
\item 
The sheaf $res(\S)$ has a natural structure of $\O_A$-module. Moreover if $\S$ is a coherent $\O_X$-module sheaf then $res(\S)$ is a coherent $\O_A$-module sheaf;
\item
There is a natural homomorphism $\S\M^{\v}\rightarrow \S(\v)$ for which we have the short exact sequence
\[
0\rightarrow\S\M^{\v+1}\rightarrow\S\M^{\v}\rightarrow\S(\v)\rightarrow 0
\]
\item
The sheaf $A_{(\v)}$ contains nilpotent elements; more precisely
\[
nil(A_{(\v)})=\{x\mid \exists n\in \N, x^n=0\}=\M/\M^{\v}
\]
\item 
There is natural isomorphisms
\[
\Q_1\cong\underline{N}^*
\]
\[
 \Qv\cong \Q_1\otimes\Q_1\otimes\cdots\otimes\Q_1 ( \v \hbox{ times})
\]
In particular $\Qv$ is a coherent $\O_A$-module sheaf (Note that $nil(\nei{A}{\v})$ is not $\O_A$-module sheaf, it is just a sheaf of $\C$-algebras) and 
\[
\Qv\cong{(\underline{N}^*)}^\v 
\]
\item 
There is natural short exact sequence
\[
0\rightarrow \Q_{\v-1}\rightarrow \nei{A}{\v}\rightarrow\nei{A}{\v-1}\rightarrow 0
\]
\item
$\nei{A}{1}=\O_A$ is the structural sheaf of $A$;
\item
There is a natural isomorphism
\[
\T/\T_A\cong (\Q_1)^*
\]
\end{itemize}
\subsection{Extension Step by Step}

Let $A\hookrightarrow X$ be an embedding of the variety $A$ in a complex manifold $X$. We want to describe the germ of this embedding with a minimal number  of datas. The first elementary data of an embedding is its normal bundle (when $A$ is not smooth the sheaf $\Q_1=\M/\M^2$ plays the role of normal bundle). The other datas of an embedding are its finite neighborhoods. The following proposition shows that the higher level neighborhoods contain the informations of lower level neighborhoods.
\begin{pro}
Any isomorphism $\nei{\phi}{\v}: A_{(\v)}\rightarrow A'_{(\v)},\v \geq 2$ induces natural
isomorphisms
\[
 A_{(\mu)}\rightarrow A'_{(\mu)},\mu\leq\v
\]
\[
\Q_1\rightarrow\Q_1'
\]
\end{pro}
We formulate our main problem in this section as follows: Let $A'$ be the image of another embedding of $A$ in a manifold $X'$.
\begin{enumerate}
\item
Given an isomorphism $\phi:\Q_1\rightarrow\Q_1'$. Under which conditions it is induced by an isomorphism $\nei{\phi}{2}:  \nei{A}{2}\rightarrow \nei{A}{2}$?
\item
Given an isomorphism $\nei{\phi}{\v}: \nei{A}{\v}\rightarrow \nei{A}{\v}, \v\geq 2$. Under which conditions it extends to $\nei{\phi}{\v+1}: \nei{A}{\v+1}\rightarrow \nei{A}{\v+1}$? 
\end{enumerate}
Note that if all such conditions are satisfied for $A$ and $A'$, we get only an isomorphism of formal neighborhoods of $A$ and $A'$.

Let $a\in A$ and $a'=\phi(a)$ be its corresponding point in $A'$. The stalk of the sheaf $\nei{A}{\v}, \v\in \tilde{N}$ at $a$ is denoted by $\nei{A}{\v}_a$. Any isomorphism
\begin{equation}
\label{111}
\nei{\phi}{\v}_{a}:\nei{A}{\v}_{a}\rightarrow \nei{A'}{\v}_{a'}
\end{equation}
determines an isomorphism  between $\nei{A}{\v}\mid_{U_a}$ and $\nei{A'}{\v}\mid_{U_{a'}}$, where $U_a$ and $U_{a'}$ are two open neighborhood of $a$ and $a'$ in $A$ and $A'$, respectively.

The following proposition gives us the local solutions of our problem:
\begin{pro}
\label{10out}
Any isomorphism $\nei{\phi}{\v}_{a}:\nei{A}{\v}_{a}\rightarrow \nei{A'}{\v}_{a'}$ is induced by an isomorphism 
\begin{equation}
\label{111}
\nei{\phi}{*}_a:\nei{A}{*}_a\rightarrow \nei{A}{*}_{a'}
\end{equation}
and hence extends to
\begin{equation}
\label{222}
\nei{\phi}{\v+1}_a:\nei{A}{\v+1}_a\rightarrow \nei{A}{\v+1}_{a'}
\end{equation}
\end{pro}
The above proposition is the third part of Proposition ~\ref{sara}.
\begin{rem}
The isomorphism $\nei{\phi}{*}:A_{*,a}\rightarrow A'_{*,a'}$ is not unique. 
\end{rem}
In the introduction of \cite{grre1} we find the following statement of 
H. Cartan: {\it la notion de faisceau s'introduit parce qu'il s'agit de passer 
de donn\'ees locales \`a l'etude de propri\'et\'es globales}.
Like many other examples in complex analysis, 
the obstructions to glue 
the local solutions lie in a first cohomology group of a sheaf over $A$. 
The precise identification of that sheaf and its first cohomology group 
is our main objective in this section.

Now, let us be given an isomorphism $\nei{\phi}{\v}:A_{(\v)}\rightarrow A'_{(\v)}$. We want to extend $\nei{\phi}{\v}$ to $\nei{\phi}{\v+1}:A_{(\v+1)}\rightarrow A'_{(\v+1)}$, i.e., to find an isomorphism $\nei{\phi}{\v+1}:A_{(\v+1)}\rightarrow A'_{(\v+1)}$ such that the following diagram is commutative:
\begin{equation}
\begin{array}{ccc}
\nei{A}{\v+1}& \stackrel{\nei{\phi}{\v+1}}{\rightarrow}  & \nei{A}{\v+1} \\
\downarrow &  & \downarrow \\
\nei{A}{\v} & \stackrel{\nei{\phi}{\v}}{\rightarrow}  & \nei{A}{\v} \\
\end{array}
\end{equation}
Proposition ~\ref{10out} gives us the local solutions 

\begin{equation}
\label{10out1}
\begin{array}{ccc}
\nei{A}{\v+1}_a& \stackrel{\nei{\phi}{\v+1}_a}{\rightarrow}  & \nei{A'}{\v+1}_{a'} \\
\downarrow &  & \downarrow \\
\nei{A}{\v}_{a} & \stackrel{\nei{\phi}{\v}_a}{\rightarrow}  & \nei{A'}{\v}_{a'} 
\end{array}
\end{equation}
where $\nei{A}{\v}_{a}$ is the stalk of the sheaf $\nei{A}{\v}$ over the point $a$. Now, cover $A$ with small open sets for which we have the diagrams of the type
~(\ref{10out1}). Combining two diagrams in the intersection of 
neighborhoods of the points $a$ and $b$ we get:

\begin{equation}
\label{10out2}
\begin{array}{ccc}
\nei{A}{\v+1}_{a,b}& \stackrel{\nei{\phi}{\v+1}_{a,b}}{\rightarrow}  & \nei{A'}{\v+1}_{a',b'} \\
\downarrow &  & \downarrow \\
\nei{A}{\v}_{a,b} & \stackrel{id}{\rightarrow}  & \nei{A'}{\v}_{a',b'} 
\end{array}
\end{equation}
where 
\begin{equation}
\label{tran}
\nei{\phi}{\v+1}_{a,b}=\nei{\phi}{\v+1}_{a}\circ\nei{\phi}{\v+1}_{b}^{-1}
\end{equation}
\begin{rem}\rm
Note that we have used the notation $\nei{\phi}{\v+1}_{a,b}$ instead of
$\nei{\phi}{\v+1}\mid_{U_a\cap U_b}$, $\nei{\phi}{\v+1}_{a}$ instead of
$\nei{\phi}{\v+1}\mid_{U_a}$ and so on.
\end{rem}
The above transition elements are obstruction to our extension problem. Now it is natural to define the following sheaf:
\begin{defi}\rm
$Aut(\v)$ is the sheaf of isomorphisms 
$\nei{\phi}{\v+1}:\nei{A}{\v+1}\rightarrow 
\nei{A}{\v+1}$ inducing the identity in $\nei{A}{\v}$, i.e. the following diagram is commutative
\begin{equation}
\label{10out2}
\begin{array}{ccc}
\nei{A}{\v+1} & \stackrel{\nei{\phi}{\v+1}}{\rightarrow}  & \nei{A}{\v+1} \\
\downarrow &  & \downarrow \\
\nei{A}{\v} & \stackrel{id}{\rightarrow}  & \nei{A}{\v} 
\end{array}
\end{equation}

\end{defi}
Later in Proposition ~\ref{dolphin} we will see that $Aut(\v)$ is a sheaf of Abelian groups.

Now it is easy to see that the data in ~(\ref{tran}) form an element of 
\[
H^1(A, Aut(\v))
\]
The elements of $H^1(A, Aut(\v))$ are obstruction to the extension problem.

It is clear that the case $\v=1$ needs an special treatment. $\nei{A}{1}$ is the structural sheaf of $A$ and  the condition 
$H^1(A, Aut(1))=0$ means that any two embedding of $A$ have the same 2-neighborhood and in particular have isomorphic $\nb$'s !. Therefore, the definition of $Aut(1)$ is not useful. We modify this  definition as follows:
\begin{defi}\rm
$Aut(1)$ is the sheaf of isomorphisms $\nei{\phi}{2}:\nei{A}{2}\rightarrow 
\nei{A}{2}$ inducing the identity on $\nb$ and for which the 
 following diagram is commutative
\begin{equation}
\begin{array}{ccc}
\nei{A}{2} & \stackrel{\nei{\phi}{2}}{\rightarrow}  & \nei{A}{2} \\
\downarrow &  & \downarrow \\
\nei{A}{1} & \stackrel{id}{\rightarrow}  & \nei{A}{1} 
\end{array}
\end{equation}
\end{defi} 

\begin{pro}
If $H^1(A, Aut(\v))=0$ then any isomorphism 
\begin{enumerate}
\item
$\nei{\phi}{\v}:A_{(\v)}\rightarrow A'_{(\v)}$ if  $\v > 1$
\item
$\nei{\phi}{}:\Q_1\rightarrow \Q_1'$ if $\v=1$
\end{enumerate}
  extends to an isomorphism $\phi_{(\v+1)}:A_{(\v+1)}\rightarrow A'_{(\v+1)}$. 
 \end{pro}

\pr
The obstruction to the above extension is obtained by diagram ~(\ref{10out2}) and so is an element of $H^1(A, Aut(\v))$.

Now we have to identify $Aut(\v)$ and especially we have to verify when 
$H^1(A, Aut(\v))=0$ is satisfied. 
\\
\begin{pro}
\label{dolphin}
For $\v\geq 2$ we have 
\[
Aut(\v)\equiv \T(\v)(\stackrel{def}{=} \T \otimes_{\O_A} \Q_{\v}) 
\]
where $\T$ is the sheaf of holomorphic vector fields in $X$ (sections of the tangent bundle of $X$); for the case $\v=1$ we have
\[
Aut(1)\equiv \T_A(1)
\]
where $\T_A$ is the sheaf of holomorphic vector fields in $X$ tangent to $A$.  
\end{pro}

\pr
Let us introduce the function which will be our candidate for the desired isomorphisms. First consider the case $\v\geq 2$. 
\[
*: \T(\v)\rightarrow Aut(\v)
\]
For any $\psi\in \T(\v)$ define
\[
\beta,\beta':\nei{A}{\v+1}\rightarrow\nei{A}{\v+1}
\]
\[
\beta(f)=f+\psi . df
\]
\[
\beta'(f)=f-\psi .df
\]
we have
\[
\beta\circ\beta' (f)=f-\psi.df+\psi d(f-\psi.df)=f-\psi.d(\psi .df)=f
\hbox{ mod } \M^{2\v-1}
\]
We have  $2\v-1\geq \v+1$ and so
\begin{equation}
\label{3nov00}
\beta\circ\beta' (f)=f \hbox{ mod } \M^{\v+1}
\end{equation}
In other words  $\beta'$ is the inverse function of $\beta$. We define
\[
*(\psi)=\beta
\]
Now it is enough to prove that $*$ is the desired isomorphism. Since $X$ is nonsingular $*$ is injective. Let $\beta\in Aut(\v)$. We write
\[
\beta(f)-f=\psi'(f)
\]
$\psi'(f)=0$ mod $\M^\v$ and so $\psi'\in Hom(\nei{A}{\v+1}, \nbv)$. Composing with $\nei{A}{*}\rightarrow \nei{A}{\v+1}$ and without change in notations we can assume
\[
\psi'\in Hom(\nei{A}{*}, \nbv)
\]
Let $z_1,z_2,\ldots,z_n$ be local coordinates. Define
\[
\psi(dz_i)=\psi'(z_i)
\]
Then $\psi\in \T(\v)$ and the mapping $\beta\rightarrow \psi$ is the inverse of $*$.

The case $\v=1$ is the same as previous one. We need to substitute $\T_A$ for $\T$  to get the congruency ~(\ref{3nov00}).\qed

How can we calculate the cohomology groups $H^1(A, \T(\v))$? To do this, we break $\T(\v)$ into two other simple sheaves as follows:
\\
There is a natural short exact sequence 
\[
0\rightarrow \T_A \rightarrow \T\rightarrow \Q_1^*\rightarrow 0
\]
By tensorial multiplication over $\O_A$ with $\Q_\v$, we have
\[
0\rightarrow \T_A(\v) \rightarrow \T(\v)\rightarrow \Q_{\v-1}\rightarrow 0
\]
This gives us the long exact sequence
\[
\ldots\rightarrow H^1(A, \T_A(\v))\rightarrow H^1(A, \T(\v))\rightarrow H^1(A, \Q_{\v-1})\rightarrow\ldots
\]
We summarize the above arguments in the following proposition:
\begin{pro}
\label{imaginacion}
If $H^1(A, \T_A(\v))=0$ and $H^1(A, \Q_{\v-1})=0$ then $ H^1(A, \T(\v))=0$ and so  any isomorphism 
\[
\phi_{(\v)}:A_{(\v)}\rightarrow A'_{(\v)}  \hbox{ if }  \v > 1
\]
\[
\phi:\Q_1\rightarrow \Q_1' \hbox{ if } \v=1
\]
extends to an isomorphism $\phi_{(\v+1)}:A_{(\v+1)}\rightarrow A'_{(\v+1)}$.
\end{pro}
\subsection{Positive and Negative Line Bundles and Vanishing Theorems }
The aim of this section is to introduce vanishing theorems in complex analysis and algebraic geometry. Let us start with the definition of negative vector bundle.
\begin{defi}
The vector bundle $V\rightarrow A$ over a complex manifold $A$
 is called negative (in the sense of Grauert) if its zero section is an exceptional set in $V$. Naturally $V\rightarrow A$ is called positive if its dual is negative.
\end{defi}
There is another definition in algebraic geometry for a positive line bundle as follows: The line bundle $L\rightarrow A$ over a complex manifold is called positive (in the sense of Kodaira) if its Chern class $c(L)\in H^2(A,\C)$ is represented by a positive $(1,1)$-form
\[
\omega=\sum g_{ij}dz_i\wedge d\overline{z_j}
\]
For more information about this definition of positive line bundles the reader is referred to \cite{gri}. 
\begin{theo}
The line bundle $L\rightarrow A$ is positive in the sense of Kodaira if and only if it is positive in the sense of Grauert.
\end{theo}
This is Satz 1 p. 341 of \cite{gra}. The various definitions of positive line bundles coincide, however, for vector bundles whose fibers have dimension greater than one these definitions are not equivalent ( see \cite{gri1}, \cite{gri2} and \cite{ume}).

Let $A$ be a compact exceptional submanifold of a manifold $X$. This situation can be obtained in negative vector bundles.  
\begin{theo}
\label{zan}
(Grauert, \cite{gra}, Hilfssatz 1, p. 344)
Let  $\S$ be a coherent analytic sheaf on $X$.
There exists a positive integer $\v_0$ such that 
\[
H^{\mu}(A, \S(\v))=0, \ \ \mu\geq 1, \ \v\geq \v_0
\]
\end{theo} 
Theorem ~\ref{26jan01} is the key of the proof of above theorem. When $X$ is negative line bundle over $A$, this theorem is exactly Kodaira's vanishing theorem. In the case where $A$ is a Riemann surface in a two dimensional complex manifold this theorem is already proved in Section ~\ref{raquel} using the Serre Duality. In many situation like this one, we are interested in the nature of the minimum number $\v_0$. The complete proof of Theorem ~\ref{zan}.
can be found in \cite{nar}. 

\begin{theo}
\label{25mar01}
(Grauert \cite{gra},Satz 2, p. 357) There exists a positive integer $k_0$ such that 
\[
H^{\mu}(U, \S\M^\v)=0, \ \ \mu\geq 1, \ \v\geq k_0
\]
where $U$ is a strongly pseudoconvex neighborhood of $A$ in $X$.
\end{theo}
\begin{rem}\rm
\label{bakht}
The following discussion shows that $\v_0\geq k_0$ can be supposed:

 Consider the short exact sequence
\[
0\rightarrow\S\M^{\v+1}\rightarrow\S\M^{\v}\rightarrow\S(\v)\rightarrow 0
\]
For $\v\geq\v_0$ the map $H^{\mu}(U, \S\M^{\v+1})\rightarrow H^{\mu}(U, \S\M^\v)$ is surjective and so for any $k\geq \v$ the map $H^{\mu}(U, \S\M^{k})\rightarrow H^{\mu}(U, \S\M^\v)$ is surjective. For a large $k$ we have $H^{\mu}(U, \S\M^{k})=0$ and therefore $H^{\mu}(U, \S\M^{\v})=0$.
\end{rem}

In what follows we will write $H^{\mu}(A, \S\M^{k})$ instead of $H^{\mu}(U, \S\M^{k})$. Here we consider sheaf theory restriction of 
$\S\M^{k}$ to $A$. 
\subsection{Grauert Theorem}
Now we are in a position to state Grauert theorem about rigidity of negatively embedded varieties. Let $A\hookrightarrow X$ and $A'\hookrightarrow X'$ be two compact submanifold of codimension one ($A\cong A'$). Let also $N$ be the normal bundle of $A$ in $X$, ${\cal D}$ its zero section and $T_A$ the tangent bundle of $A$.
\begin{theo}(Grauert \cite{gra} Satz 7 p. 363) 
\label{grauert}
Let $\nei{\phi}{l}:\nei{A}{l}\rightarrow\nei{A'}{l}$ be an isomorphism and $H^1(A, T_A\otimes N^\v)=0, H^1(A, N^{\v-1})=0\ \v\geq l$. Then $\phi$ extends to a biholomorphism of neighborhoods 
$\nei{A}{*}\rightarrow\nei{A'}{*}$.

In particular if $H^1(A, T_A\otimes N^\v)=0, H^1(A, N^\v)=0,\ \v\geq 1$ then there exists a biholomorphism between a neighborhood of $A$ in $X$ and a neighborhood of ${\cal D}$ in $N$.
\end{theo}
\pr
By the hypotheses and Proposition ~\ref{imaginacion} we can get a formal isomorphism of $(X,A)$ and $(X',A')$. Now by formal principle for negatively embedded varieties Theorem ~\ref{25jan01} we can find the desired biholomorphism. \\ 

\qed

Note that by Kodaira vanishing theorem (Theorem ~\ref{zan}) there exists a $\v_0$ such that 
$H^1(A, T_A\otimes N^\v)=0, H^1(A, N^{\v-1})=0\ \v\geq \v_0$. Roughly speaking, the germ of a negatively embedded compact manifold of codimension one is determined by a $\v$-neighborhood for $\v$ enough big.
\subsection{The Case of a Riemann Surface Embedded in a 2-Manifold}
\label{raquel}
In this section we want to substitute the vanishing conditions
$H^1(A, \T_A(\v))=0$ and $H^1(A, \Q_{\v-1})=0$ for some numerical ones. The Serre duality will be used for this purpose.
 \begin{theo}
(Serre Duality) Let $A$ be a complex manifold of complex dimension $n$ and $V$ a holomorphic vector bundle over $A$. Then there exists a natural $\C$-isomorphism 
\[
H^q(A,\underline{\Omega^p\otimes V})\cong (H^{n-q}(A,\underline{\Omega^{n-p}\otimes V^*}))^*
\]
\end{theo}
From now on we do not use the line under bundles (it denotes the sheaf of sections), for instance instead of $H^1(A, \underline{\Omega})$ we write $H^1(A, \Omega)$. Let $A$ be a Riemann surface. Putting $p=0,\ q=0$ we have 
\[
H^1(A,V)\cong (\Gamma (A,\Omega\otimes V^*))^*
\]

Now
\[
H^1(A, \Q_{\v-1})=H^1(A, (N^*)^{v-1}))=(\Gamma (A,\Omega\otimes N^{v-1}))^*
\]
$\Omega\otimes N^{v-1}$ has not any global holomorphic section if 
\begin{equation}
\label{5nov2000}
c(\Omega\otimes N^{v-1})=2g-2+(\v-1)A.A<0 
\end{equation}
In the same way 
\[
H^1(A, \T_A(\v))=(\Gamma (A,\Omega\otimes (TA)^*\otimes N^\v))^*=(\Gamma (A,\Omega\otimes \Omega\otimes N^\v))^*=0
\]
 if
\begin{equation}
\label{5nov00}
c(\Omega\otimes \Omega \otimes N^\v)=2(2g-2)+\v A.A<0
\end{equation}
Finally we conclude that 
\begin{theo}
\label{26jan27mai}
Let $A$ be a Riemann surface of genus $g$ embedded in a two dimensional manifold $X$. Suppose that 
\begin{itemize}
\item
$A.A\leq 0$ if $g=0$;
\item
$A.A<2(2-2g)$ if $g\geq 1$ 
\end{itemize}
Then the embedding $A\hookrightarrow X$ is formally equivalent with $A'\hookrightarrow X'$, where the normal bundle of $A'$ in $X'$ equals 
the normal bundle of $A$ in $X$. In the case where $A.A<0$, these two germs are holomorphically equivalent.  
\end{theo}

\pr 
Since the normal bundle of $A'$ in $X'$ equals 
the normal bundle of $A$ in $X$, there exists an isomorphism $\nei{\phi}{}:\Q_1\rightarrow \Q_1'$. To extend this isomorphism to a formal isomorphism of the neighborhoods of $A$ and $A'$ in $X$ and $X'$, respectively, we must have the inequalities ~(\ref{5nov00}) for all $\v\geq1$ and ~(\ref{5nov2000}) for all $\v>1$ satisfied.
This implies exactly $A.A\leq 0$ if $g=0$ and $A.A<2(2-2g)$ if $g\geq 1$. If $A.A<0$ then $A$ is exceptional and so  it satisfies the formal principle. The proposition is proved.
\subsection{Embedding Dimension}
Any germ of a singularity $(X,x)$ is an analytic subspace of $\C^n$, for some $n$. The smallest integer $n$ with this property is called the embedding dimension of $X$ at $x$ and 
is denoted by $emb_x X$. Recall that $\O_{X,x}$ is the ring of germs of holomorphic functions 
in a neighborhood of $x$ and $\M_x$ is the maximal ideal of the local ring $\O_{X,x}$. The $\O_{X,x}$-module $\frac{\M_x}{\M_x^2}$ is a finite dimensional vector space and we have: 
\begin{pro}(\cite{grre1} p. 115)
\label{goftam}
$emb_x{X}=dim_\C \frac{\M_x}{\M_x^2}$ for all points $x\in X$.
\end{pro}
Let $A\hookrightarrow X$ be a germ of negative embbeding and $(Y,y)$ be the singularity obtained after blowing down $A$. We are interested in the embedding dimension of $(Y,y)$ in terms of the information which we know about the embedding $A\hookrightarrow X$.

Suppose that $A$ and $X$ are complex manifolds and $N$ is the normal bundle of $A$ in $X$.
\begin{pro}
We have
\[
emb_y(Y)\leq dim_\C H^0(A,N^*)
\]
The equality happens if $H^1(A, N^{*k})=0$ for all $k\geq 1$ 
\end{pro}
\pr
There is a natural isomorphism between $\frac{\M_y}{\M_y^2}$ and $\frac{H^0(A, \M)}{H^0(A, \M^2)}$. Now the first statement can be derived from the long exact sequence 
\[
0\rightarrow H^0(A, \M^2)\rightarrow H^0(A, \M)\rightarrow H^0(A, \frac{\M}{\M^2})
\rightarrow 
\]
\[
H^1(A, \M^2)\rightarrow H^1(A, \M)\rightarrow H^1(A, \frac{\M}{\M^2})\rightarrow \ldots 
\]
of the short excact sequence $0\rightarrow \M^2 \rightarrow \M \rightarrow \frac{\M}{\M^2}\rightarrow 0$. Recall that the sheaf of holomorphic sections of $N^*$ is isomorphic to $\frac{\M}{\M^2}$.
\\
By the discussion after Theorem ~\ref{25mar01} and by the hypothesis we conclude that
\[
H^1(A,\M^k)=0,\ k\geq 1
\]
In particular $H^1(A,\M^2)=0$ implies that 
\[
\frac{H^0(A, \M)}{H^0(A, \M^2)}=H^0(A, \frac{\M}{\M^2})=H^0(A, N^*)
\]
\qed

\section{ Foliated Neighborhoods}
\def\Z{{\Bbb Z}}
\def\p1{{\Bbb P}_1}
\def\pn{{\Bbb P}_n}
\def\F{{\cal F}}
\def\Sf{\F(X,A)}
Let $A$ be a Riemann surface embedded in a two dimensional manifold $X$ with negative self-intersection. By $(X,A)$ we mean a small neighborhood of $A$ in $X$ (the germ of embedding).
In this section we consider the germs of holomorphic foliations in a 
neighborhood of $A$. Two foliations $\F$ and $\F'$ in $(X,A)$ are called equivalent if there exists a biholomorphism 
\[
\psi: (X,A)\rightarrow (X,A)
\]
such that 
\begin{enumerate}
\item
$\psi\mid_{A}$ is identity;
\item
$\psi^{-1}(\F')=\F$.
\end{enumerate}
We are interested in the space $\Sf$ of equivalence classes. Natural questions in this direction arise:
\begin{enumerate}
\item
Is $\Sf$ finite dimensional?
\item
Does $\Sf$ has a natural structure of complex space?
\item
When $\Sf$ is a discrete set? 
\end{enumerate}
 To ask these questions we start with the most simple foliations, namely, foliations without singularity and  transverse to $A$. We remember that by Theorem ~\ref{26jan27mai} if $A.A<min\{0,2(2-2g)\}$ then the germ of $A$ in $X$ is biholomorphically equivalent with the germ of zero section of its normal bundle.
\subsection{Transversal Foliations}
\def\L{{\cal L}}
In this section we are concerned with germs of transverse holomorphic 
foliations in $(X,A)$, i.e. the foliations with no singularity and with 
leaves transverse to $A$. Let us introduce some examples in the case $A=\p1$.
\begin{figure}[t]
\begin{center}
\includegraphics{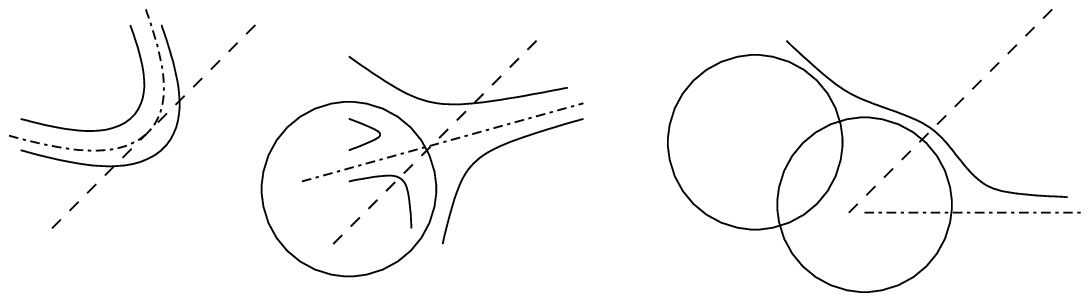}
\caption{}
\label{b}
\end{center}
\end{figure}
\begin{exa}\rm
By successive blow-ups at the origin of $\C^2$, we can get a $A\cong \p1$ embedded in a two dimensional manifold and with $A.A=-n$. A neighborhood of $A$ is covered by coordinate systems
$(u,y)=(\frac{X}{Y}, Y)$
and $(x,t)=(\frac{X^n}{Y^{n-1}},\frac{Y}{X})$, where $X$ and $Y$ are the pullback of a coordinates 
system at the origin of $\C^2$. The change of coordinates is given by
\[
(x,t)\rightarrow (\frac{1}{t},xt^n)=(u,y)
\]
 In this example we have a germ of transverse  holomorphic foliation $\F$ given by the 1-form
\[
\omega=XdY-YdX=(xt^{n-1})^2dt=-y^2du
\]
It is easy to check that
\[
zer(\omega)=2.A+2(n-1)L
\]
\[
zer(Y)=1.A+n.L, \  zer(X)=1.A+(n-1)L+L'
\]
where $zer()$ means the zero divisor and $L$ (resp. $L'$) is the leaf of $\F$ given by $t=0$ (resp. $u=0$ ) in the 
coordinates $(x,t)$ (resp. $(u,y)$); it is the pullback of $X$-axis (resp. $Y$-axis). 
\end{exa}
The above example contains the basic idea of the proof of the following theorem.

\begin{theo}
\label{27jan01}
Let $A$ be a Riemann surface of genus $g$ embedded in a 
manifold $X$ of dimension two with $A. A<min\{2-2g,0\}$. 
The germs of any two holomorphic  transverse foliations  are equivalent.    
\end{theo}
The group of line bundles in a complex manifold $M$ is called the Picard group and is denoted by $Pic(M)$. 
Let us first state the main lemma we need in the proof of the above theorem:
\begin{lem}
\label{nashenas}
Let $A$ be a Riemann surface of genus $g$ embedded in a 
manifold $X$ of dimension two with $A. A<min\{2-2g,0\}$. The restriction map
$r:Pic(X)\rightarrow Pic(A)$ is injective. Moreover, if $X$ has transverse foliation to $A$ then $r$ is an isomorphism.
\end{lem}
\pr
By Serre duality we can see that the condition $A. A<min\{2-2g,0\}$ implies that
\[
H^\mu(A, \M/\M^2)=0,\ \mu\geq 1
\]
By Remark ~\ref{bakht} we conclude that
\[
H^\mu(U, \M)=0,\ \mu\geq 1
\]
where $U$ is a strongly pseudoconvex neighborhood of $A$ in $X$.
The diagram
\begin{equation}
\begin{array}{ccccccccc}
 & & & & 0  & & & &  \\
& & & & \downarrow & & & &  \\
& & & & \M & & & &  \\
& & & & \downarrow & & & &  \\
0 &\rightarrow & \Z & \rightarrow&\O_X &\rightarrow &\O_X^* &\rightarrow & 0  \\
& &\downarrow & &\downarrow & &\downarrow & &  \\
0 &\rightarrow & \Z & \rightarrow&\O_A &\rightarrow &\O_A^* &\rightarrow & 0  \\
& & & & \downarrow & & & &  \\
 & & & & 0  & & & &  \\
\end{array}
\end{equation}
gives us 
\begin{equation}
\begin{array}{ccccccccc}

& &  H^1(U, \M)=0& & & & & &  \\
& & \downarrow& &  & & & &  \\
H^1(U,\Z) &\rightarrow & H^1(U, \O_X) & \rightarrow& H^1(U, \O_X^*) &\rightarrow & H^2(U,\Z) &\rightarrow & \cdots  \\
\downarrow & &\downarrow & &\downarrow & &\downarrow & &  \\
 H^1(A, \Z) &\rightarrow & H^1(A, \O_A) & \rightarrow& H^1(A, \O_A^*) &\rightarrow &H^2(A,\Z) &\rightarrow & \cdots  \\
& & \downarrow& &  & & & &  \\
 & &  H^2(U, \M)=0& & & & & &  \\
\end{array}
\end{equation}
By considering a small neighborhood $U$, if necessary, we can assume that $A$ and $U$ have the same topology and so the first and forth column functions are isomorphism. The second column is also an isomorphism.  Now a simple argument
shows that the third column is injective.
\\
If $(X,A)$ has transverse foliation then $r$ is surjective and so it is an isomorphism.
\qed

{\bf Proof of Theorem ~\ref{27jan01}:}
Let $\F$ be the germ of a transverse foliation in $(X,A)$ and $N$ the normal bundle of $A$ in $X$. The normal bundle $N$ of $A$ in $X$
 has a meromorphic global section namely $s$. Let
\[
div(s)=\sum n_ip_i,\ p_i\in A,\ n_i\in\Z
\]
We define the divisor $D$ in $X$ as follows:  
\[
D=A-\sum n_i\L_{p_i}
\]
where $\L_{p_i}$  is the leaf of $\F$  through $p_i$.  
The line bundle $L_D$ associated to $D$  restricted to $A$ is trivial line bundle, and so by Lemma ~\ref{nashenas}, $L_D$ is trivial or equivalently 
there exists a meromorphic function $g$ on $(X,A)$ with $div(g)=D$.

Let $\tilde{f}$ be an arbitrary meromorphic function on $A$ and $f$ its extension along the foliation . Define the 1-form 
\[
\omega=gdf
\]
The 1-form $\omega$ has the following properties
\begin{enumerate}
\item
$\omega$ induces the foliation $\F$;
\item
The divisor of $\omega$ is $A+K$, where $K$ is $\F$-invariant and its restriction to $A$ depends only on $\tilde{f}$ and the meromorphic section $s$. 
\end{enumerate}
Let $\F'$ be another transverse foliation in $(X,A)$. In the same way we can construct the 1-form $\omega'$ for $\F$. By taking normalizing charts at $a\in A$ for $\F$ and $\F'$, we can verify easily that there exists a unique biholomorphism 
\[
\psi_a:(X,A,a)\rightarrow (X,A,a)
\]
inducing identity on $A$ and with the property $\psi^{-1}(\omega')=\omega$. Now the uniqueness of these local biholomorphisms implies the existence of a biholomorphic equivalence between $\F$ and $\F'$.
\qed
\subsection{Foliations with Tangencies}
\begin{figure}[t]
\begin{center}
\includegraphics{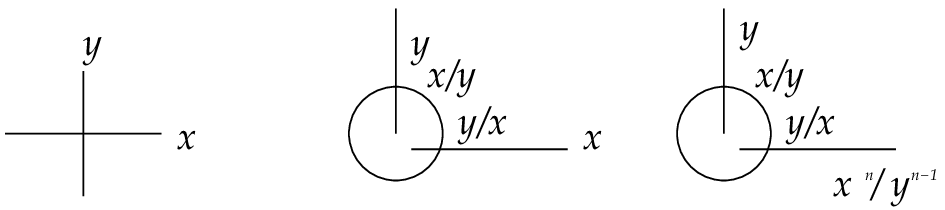}
\caption{}
\label{b}
\end{center}
\end{figure}
We start this section with M. Suzuki's examples (\cite{suz} and \cite{suz1}). Consider The germ of holomorphic foliations given by 
\[
\F(\omega): \omega=(y^3+y^2-xy)dx-(2xy^2+xy-x^2)dy=0
\]
\[
\F(\eta): \eta=(2y^2+x^3)dx-2xydy
\]
In both cases after blowing up at $0$ we get two non singular foliations around $\p1$, the divisor of blow up, and with the following property: The foliation is transverse to $\p1$ in all points except one point and in this point it has a tangency of order two with $\p1$.
It is shown in \cite{cermat} that these two foliations are topologically equivalent but not holomorphically.
\subsection{Formal Isomorphism of foliations}
The aim of this section is to extend the methods used in the section ~\ref{10feb01} to the case where we have foliated neighborhoods. We define the formal isomorphism between two foliated neighborhoods and then we identify the obstructions for the existence of such formal isomorphism. We prove that if the foliation $\F$ has not singularities on $A$ and 
$A.A<min\{0,2-2g-tang(\F,A)\}$ then any other holomorphic foilation in a neighborhood of $A$ having the same local analytic structure of $\F$, is formally isomorphism with $\F$. Here $tang(\F,A)$ is the number of tangency points between $A$ and $\F$ counting with
multiplicity.
 In another words the formal moduli space of foliations with the local structure of $\F$ contains only one point.

Let $\F$ and $\F'$ be the germs of two foliations in $(X,A)$ which are locally biholomorphic, i.e. for any point $a\in A$ there exists a biholomorphism 
\[
\phi_a :(X,A,a)\rightarrow (X,A,a)
\]
sending the foliation $\F$ to $\F'$. Roughly speaking, $\F$ and $\F'$ has the same local analytic structure around $A$.
\begin{defi}\rm
Let $\v$ be a natural number. 
We say thet the ismorphism  
\[
\nei{\phi}{\v}: \nei{A}{\v}\rightarrow \nei{A}{\v}
\]
is $\v$-isomorphism between $\F$ and  $\F'$ if for every point $a\in A$ 
there exists a local biholomorphism 
\[
\phi_a :(X,A,a)\rightarrow (X,A,a)
\]
which induces $\nei{\phi}{\v}$ and sends  $\F$ to $\F'$. We also say that $\nei{\phi}{\v}$ sends $\F$ to $\F'$.
\\
The isomorphism of formal neighborhoods 
\[
\nei{\phi}{\infty}:\nei{A}{\infty}\rightarrow \nei{A}{\infty}
\]
is a formal biholomorphism between $\F$ and $\F'$ if for every natural number $\v$ the $\v$-isomorphism $\nei{A}{\v}$ induced by $\nei{\phi}{\infty}$ sends $\F$ to $\F'$.
\end{defi}  
Now we are going to identify the obstructions for the existence of formal biholomorphism between two foliations.

Let us be given an $\v$-isomorphism 
$\nei{\phi}{\v}:\nei{A}{\v}\rightarrow \nei{A}{\v}$ between the foliations $\F$ and $\F'$. 
We want to extend $\nei{\phi}{\v}$ to $\nei{\phi}{\v+1}:A_{(\v+1)}\rightarrow A'_{(\v+1)}$, i.e. to find a $(\v+1)$-isomorphism $\nei{\phi}{\v+1}:A_{(\v+1)}\rightarrow A'_{(\v+1)}$ between $\F$ and $\F'$
 such that the following diagram is commutative:
\begin{equation}
\begin{array}{ccc}
\nei{A}{\v+1}& \stackrel{\nei{\phi}{\v+1}}{\rightarrow}  & \nei{A}{\v+1} \\
\downarrow &  & \downarrow \\
\nei{A}{\v} & \stackrel{\nei{\phi}{\v}}{\rightarrow}  & \nei{A}{\v} \\
\end{array}
\end{equation}
$\F$ and $\F'$ have the same local analytic structure. Therefore we have the local solutions of our problem.
\begin{equation}
\label{10fe1}
\begin{array}{ccc}
\nei{A}{\v+1}_a& \stackrel{\nei{\phi}{\v+1}_a}{\rightarrow}  & \nei{A}{\v+1}_{a} \\
\downarrow &  & \downarrow \\
\nei{A}{\v}_{a} & \stackrel{\nei{\phi}{\v}_a}{\rightarrow}  & \nei{A}{\v}_{a} 
\end{array}
\end{equation}
where $\nei{A}{\v}_{a}$ is the stalk of the sheaf $\nei{A}{\v}$ over the point $a$. If $\v=1$ we can furtheremore assume that $\nei{\phi}{\v+1}_a$ is an identitiy on $\M/\M^2$.  Now, cover $A$ with small open sets for which we have the diagrams of the type
~(\ref{10fe1}). Combining two diagrams in the intersection of 
neighborhoods of the points $a$ and $b$ we get:

\begin{equation}
\label{10f}
\begin{array}{ccc}
\nei{A}{\v+1}_{a,b}& \stackrel{\nei{\phi}{\v+1}_{a,b}}{\rightarrow}  & \nei{A}{\v+1}_{a,b} \\
\downarrow &  & \downarrow \\
\nei{A}{\v}_{a,b} & \stackrel{id}{\rightarrow}  & \nei{A}{\v}_{a,b} 
\end{array}
\end{equation}
where 
\begin{equation}
\label{tamam}
\nei{\phi}{\v+1}_{a,b}=\nei{\phi}{\v+1}_{a}\circ\nei{\phi}{\v+1}_{b}^{-1}
\end{equation}
sends the foliation $\F$ to itself.
\begin{rem}\rm
Note that we have used the notation $\nei{\phi}{\v+1}_{a,b}$ instead of
$\nei{\phi}{\v+1}\mid_{U_a\cap U_b}$, $\nei{\phi}{\v+1}_{a}$ instead of
$\nei{\phi}{\v+1}\mid_{U_a}$ and so on.
\end{rem}
The above transition elements are obstruction to our extension problem. Now it is natural to define the following sheaf:
\begin{defi}\rm
$Aut(\v,\F)$ is the sheaf of $(\v+1)$-isomorphisms 
$\nei{\phi}{\v+1}:\nei{A}{\v+1}\rightarrow 
\nei{A}{\v+1}$ which sends $\F$ to itself and induces the identity in $\nei{A}{\v}$, i.e. the following diagram is commutative
\begin{equation}
\begin{array}{ccc}
\nei{A}{\v+1} & \stackrel{\nei{\phi}{\v+1}}{\rightarrow}  & \nei{A}{\v+1} \\
\downarrow &  & \downarrow \\
\nei{A}{\v} & \stackrel{id}{\rightarrow}  & \nei{A}{\v} 
\end{array}
\end{equation}
in the case $\v=1$ we assume furthermore that $\nei{\phi}{\v+1}$ is an identity on $\M/\M^2$. 
\end{defi}
Now it is easy to see that the data in ~(\ref{tamam}) form an element of 
\[
H^1(A, Aut(\v,\F))
\]
The elements of $H^1(A, Aut(\v,\F))$ are obstruction to the extension problem. More precisely we have proved the following proposition:
\begin{pro}
If $H^1(A, Aut(\v,\F))=0$ then any $\v$-isomorphism  between the foliation $\F$ and $\F'$ extends to a $(\v+1)$-isomorphism between them. 
 \end{pro}
Now we have to identify $Aut(\v, \F)$ and especially we have to verify when 
$H^1(A, Aut(\v,\F))=0$ is satisfied. 
\begin{pro}
If $A$ is not $\F$-invariant then $Aut(1,\F)_a=0$ for all points $a$ in which $\F$ is 
transverse to $A$ and so $H^1(A, Aut(1,\F))=0$.  
\end{pro}
\pr
Let $\F$ be transverse  to $A$ at $a$. Choose a  coordinates system $(x,y)$ around $a$ such that $\F$ in this coordinates system is given by $x=const.$. Now it is easy to see that every biholomorphism $(\C^2,0)\rightarrow(\C^2,0)$ which sends $\F$ to $\F$ and induces the 
 identity on $\M/\M^2$ has the form
\[
(x,y)\rightarrow (x,y+y^2s_2(x)+h.o.t. )
\]
and hence induces the identity in $\nei{A}{2}$.\qed

The above proposition says when $A$ is not $\F$-invariant we can always 
find a 2-isomorphism between the foliations $\F$ and $\F'$.
\begin{theo}
\label{dolphin}
Assume that $A$ is not $\F$-invariant and $\F$ does not have singularities on $A$. 
For $\v\geq 2$ we have 
\[
Aut(\v,\F)\cong \T_\F(\v) 
\]
where $\T_\F$ is the sheaf of holomorphic vector fields in $X$ inducing the foliation $\F$.
\end{theo}
\pr
Recall that  $\T_\F(\v)=\T_\F \otimes_{\O_A} \Q_{\v}=\T_\F.\M^\v/\T_\F.\M^{\v+1}$. 
Let us introduce our candidate for the isomorphim:
\[
*: \T_\F(\v)\rightarrow Aut(\v,\F)
\]
The operator $*$ associate to every holomorphic vector field $X\in \T_\F(\v)$ the $(\v+1)$-isomorphism
\[
*(X): \nei{A}{\v+1}\rightarrow \nei{A}{\v+1}
\]
\[
f\rightarrow f+df.X
\]
Since $X$ has zero of order $\v$ in $A$, $*(X)$ induces identity in $\nei{A}{\v}$. We must prove  that $*(X)$ sends $\F$ to $\F$. 

Let $X_t(x)$ be the solution of the vector field $X$ passing through $x$ in the time $t$. Since $X$ is zero in $A$, $X_1=X_t\mid_{t=1}$ is well-defined in a smaller neighborhood around $A$. $X$ is tangent to the foliation and so $X_1$ sends $\F$ to $\F$. It is enough to prove that $X_1$ induces the map $*(X)$ in $\nei{A}{\v+1}$. We have
\[
X_t^*f=f\circ X_t=f+tdf(X)+\sum_{i\geq 2} \frac{\partial^i(f\circ X_t)}{\partial t^i}\mid_{t=0}t^i
\]
 Since
\[
\frac{\partial^2(f\circ X_t)}{\partial t^2}=
((d^2f\circ X_t).(X\circ X_t)).(X\circ X_t)+(df\circ X_t).((dX\circ X_t).(X\circ X_t))       
\]
$X$ has zero of order $\v$ along $A$ and $v\geq 2$, we conclude that
\[
\frac{\partial^i(f\circ X_t)}{\partial t^i}\mid_{t=0}=0 \hbox{ mod } \M^{\v+1} 
\]
or equivalentely
\[
X_1^*f=f+df(X)\hbox{ mod } \M^{\v+1}
\]
$*$ is trivially injective. Let us now  prove that $*$ is surjective.
\\
Let $\beta\in Aut(\F,\v)$ and 
\[
h:(x,y)\rightarrow (x,y)+(f,g)
\]
be an isomorphism in a coordinate system $(x,y)$ around a point $a\in A$ which extends $\beta$ and sends 
$\F$ to $\F$.We have  $f,g\in\M^{\v}$. Suppose that in this coordinates system $\F$ is given by the 1-form $\omega=Pdy-Qdx=0$, where $P$ and $Q$ are relatively prime. Since $h^*(\omega)\wedge\omega=0$ we have
\begin{equation}
\label{zibayi}
P\tilde{Q}f_x+Q\tilde{Q}f_y-P\tilde{P}g_x-\tilde{P}Qg_y=0
\end{equation}
where
\[
\tilde{P}=P(x+f,y+g),\ \tilde{Q}=Q(x+f,y+g)
\]
Since $A$ is not $\F$-invariant, $y$ does not devide $Q$. Therefore considering the equality ~(\ref{zibayi}) modulo $\M^\v$ we see that
\[
Qf_y-Pg_y\hbox{ mod } \M^\v
\]
This implies that 
\[
Qf-Pg=0\hbox{ mod } \M^{\v+1}
\]
The foliation $\F$ has not singularity at $a$ and so $PQ(a)\not=0$. 
Using this fact we can fined new holomorphic functions $\tilde{f}$ and $\tilde{g}$ 
such that
\[
\tilde{f}=f,\ \tilde{g}=g \hbox{ mod } \M^{\v}
\]
\[
\omega(X)=Q\tilde{f}-P\tilde{g}=0
\]
where $X=(\tilde{f},\tilde{g})$. The vector field $X$ is the desired.
\qed 


Now suppose that $A$ is not $\F$-invariant. $\F$ is transverse to $A$ except in a finite
 number of points. These points may be  tangency points of $\F$ with $A$ or singularities 
of $\F$. Suppose that there does not exists a singularity of $\F$ on $A$.

Using Serre duality, we have 
\[
H^1(A, \T_\F(\v))=(\Gamma (A,\Omega\otimes T_\F^*\otimes N^\v))^*=0
\]
 if
\begin{equation}
\label{negaran}
c(\Omega\otimes T_\F^* \otimes N^\v)=(2g-2)-c(T_\F)+\v A.A<0
\end{equation}
We have
\begin{equation}
\label{6mar01}
c(T_\F)=A.A-tang(\F,A)
\end{equation}
where $tang(\F,A)$ is the number of tangency points of $\F$ and $A$, counting with multiplicity ( see \cite{bru}).
Now substituting ~(\ref{6mar01}) in ~(\ref{negaran}), we conclude that:
\begin{theo}
\label{26jan}
Let $A$ be a Riemann surface of genus $g$ embedded in a two dimensional manifold $X$ and $\F$ and $\F'$ be two locally biholomorphic and without singularity foliations around $A$. 
If $A.A<min\{0,2-2g-tang(\F,A)\}$ the there exists a formal isomorphism between $\F$ and $\F'$. 
\end{theo}

\end{document}